\documentclass[12pt]{amsart}

\usepackage{latexsym}
\usepackage{a4}
\usepackage{amssymb}
\usepackage{amsbsy}

\newtheorem{thm}{Theorem}[section] 
\newtheorem{pro}[thm]{Proposition} 
\newtheorem{lem}[thm]{Lemma} 
 
\newtheorem{cor}[thm]{Corollary} 
\theoremstyle{definition}
 
\newtheorem{de}[thm]{Definition}
\numberwithin{equation}{section} 
\newcommand{\ab}[1]{{\mathbf{#1}}} 
\newcommand{\ob}[1]{{\mathbb{#1}}}

\newcommand{\VecTwo}[2]{ 
   \left( 
   \begin{smallmatrix} 
      #1 \\ #2 
   \end{smallmatrix} 
   \right) 
   } 
\newcommand{\N}{\Bbb{ N}}

\newcommand{\Mod}{\mathrm{Mod}} 
\newcommand{\setsuchthat}{\,\, \pmb{|} \,\,}

\newcommand{\vb}[1]{\mathbf{#1}} 
\newcommand{\cb}[1]{#1}

\newcommand{\meet}{\wedge}
\newcommand{\join}{\vee}

\newcommand{\ACC}{{\mathrm{(ACC)}}}
\newcommand{\DCC}{{\mathrm{(DCC)}}}

\newcommand{\ltlex}{<_{\mathrm{lex}}}
\newcommand{\lelex}{\le_{\mathrm{lex}}}
\newcommand{\lea}{\le_{\mathrm{E}}}

\newcommand{\FO}{\textsf{firstOcc}\,}
\newcommand{\Tab}{\tau_{\vb{a},\vb{b},h}}
\newcommand{\ran}{\textrm{range}}

\newcommand{\lcover}{\prec}

\newcommand{\algop}[2]{\langle {#1}, {#2} \rangle}

\newcommand{\U}{\mathcal{U}}

\newcommand{\PhiX}{\boldsymbol{\varphi}}

\newcommand{\Th}{\mathrm{Th}}
\newcommand{\F}{\mathcal{F}}
\newcommand{\Var}{\mathbb{V}}
\newcommand{\VarH}{\mathbb{H}}
\newcommand{\VarS}{\mathbb{S}}
\newcommand{\VarP}{\mathbb{P}}

\title{Finitely generated equational classes}
\author{Erhard Aichinger}
\author{Peter Mayr}
\thanks{Supported by the Austrian Science Fund (FWF): P24077 and P24285}
\date{\today}

\begin{document}
\bibliographystyle{amsalpha}

\maketitle
\begin{abstract}
 Classes of algebraic structures that are defined by equational laws are called \emph{varieties}
 or \emph{equational classes}.
 A variety is finitely generated if it is defined by the
 laws that hold in some fixed finite algebra. 
 We show that
 every subvariety of a finitely generated congruence permutable
 variety is finitely generated; in fact, we prove the
 more general result that
 if a finitely generated variety has an edge term,
 then all its subvarieties are finitely generated as well. 
  This applies in particular to
 all varieties of groups, loops, quasigroups and their expansions
 (e.g., modules, rings, Lie algebras, \dots).
\end{abstract}

\section{Introduction}
 Algebraic structures (or simply \emph{algebras}) are often classified according to the equational
 laws they fulfill. For example, groups can be defined as those algebras
 with operation symbols $\cdot, \mbox{}^{-1}, 1$ that satisfy the identities
 $(x \cdot y) \cdot z \approx x \cdot (y \cdot z)$, $1 \cdot x \approx x$, and $x^{-1} \cdot x \approx 1$.
 A class of algebraic structures that
 is defined by equational laws is called an \emph{equational class}, or simply
 a \emph{variety}.
 We will consider equational classes that are generated by a single finite algebra.
 As an example, the variety generated by the symmetric group $S_3$ on $3$ letters consists of those
 groups that satisfy all equational laws that are valid in $S_3$, such as
 $x^6 \approx 1$ or $x^2 \cdot y^2 \approx y^2 \cdot x^2$. 
 Since the groups in a finitely generated variety have finite exponent, clearly not every variety
 is finitely generated. However, from the work of Oates and Powell~\cite{OP:IRIF}, we obtain the
 following surprising fact for every finitely generated variety $V$ of groups:
\begin{multline} \label{eq:prop}
 \text{All subvarieties of $V$ are finitely generated.}
\end{multline}
 Kruse and L'vov proved that
 the same property is also satisfied for finitely generated varieties of rings
 \cite{Kr:ISFR,Lv:VARI}.
 For finitely generated varieties
 of lattices,~\eqref{eq:prop} is a consequence 
 of J{\'o}nsson's Lemma~\cite[Corollary~IV.6.10]{BS:ACIU}.
 
 In the main result of this paper, Theorem~\ref{thm:accsv2},  we establish that property~\eqref{eq:prop}
 holds for a much larger class of algebraic structures that includes groups, loops, rings and lattices
 as well as their expansions.

 We continue with a discussion of those concepts from equational logic that we will use.
 First we fix one \emph{type} of algebraic structures, this means, one set of operation symbols.
   An \emph{identity} (of this type) is a formula of the form $s(x_1,\ldots, x_n) \approx t(x_1,\ldots, x_n)$, where
   $s,t$ are terms involving these operation symbols.
   The fundamental relation between \emph{algebraic structures} and \emph{identities}
   is whether the algebraic structure $\ab{A}$ \emph{satisfies} the identity $s(x_1,\ldots, x_n) \approx t(x_1,    \ldots, x_n)$, which is the case 
 if $s^{\ab{A}} (a_1,\ldots, a_n) = t^{\ab{A}} (a_1,\ldots, a_n)$ for the induced term functions
 $s^{\ab{A}},t^{\ab{A}}$ on $\ab{A}$ and for all $a_1,\ldots, a_n \in A$.
   In this case, we write $\ab{A} \models s \approx t$.
   Now for a set of identities $\Phi$, the class $\Mod (\Phi)$ is defined to be the
   class of all algebras
    such that
    $\ab{A} \models \varphi$ for all $\varphi \in \Phi$.
    It is not hard to see that $\Mod (\Phi)$ is closed under taking subalgebras,
    direct products, and homomorphic images. The fundamental HSP-Theorem of
    Birkhoff \cite[Theorem~10]{Bi:Comp} tells that every class of 
    algebras closed under taking subalgebras, direct products, and homomorphic
    images is of the form $\Mod (\Phi)$ for some (possibly infinite)
    set of identities. A class of algebraic structures 
     that is of the form $\Mod (\Phi)$ is
    called an \emph{equational class} or a \emph{variety}.
    For an algebra $\ab{A}$, we consider the smallest variety that contains $\ab{A}$.
    This variety, $\Var (\ab{A})$, is given by
     \begin{equation*}
         \Var (\ab{A}) = \{ \ab{B} \setsuchthat 
     \text{for all identities } \varphi \text{ with }
                                                \ab{A} \models \varphi, \text{ we have } \ab{B} \models
                                                 \varphi \}.
    \end{equation*}
    By Birkhoff's Theorem,
    $\Var (\ab{A}) = \VarH \VarS \VarP (\ab{A})$, where $\VarH \VarS \VarP (\ab{A})$ stands
    for the class of homomorphic images of subalgebras of direct powers of $\ab{A}$.
    We call an arbitrary variety $V$ of algebras \emph{finitely generated} if there
    is a finite algebra $\ab{A}$ with $V = \Var (\ab{A})$.
     For a finitely generated variety, we can effectively check 
     whether an identity holds in all of its algebras, since it suffices to verify that
     it holds in the finite generating algebra.
    More details on varieties and identities can be found, e. g., in \cite{BS:ACIU, Mc:AFGT}.

 In this paper we consider finitely generated varieties $V$ that satisfy property~\eqref{eq:prop},
 that is, all subvarieties of $V$ are finitely generated as well. That this condition is indeed
 nontrivial can be seen from the following examples:
 Oates MacDonald and Vaughan-Lee showed that the three-element Murskii groupoid $\ab{M}$
 with zero generates a variety with subvarieties that are not finitely
 generated~\cite[Corollary~4.2]{OV:VTMO}.
 Lee proved that property~\eqref{eq:prop} does not hold for the variety generated by the
 five-element Brandt semigroup~\cite[Proposition 6.7]{Le:SVGF}.

   However,~\eqref{eq:prop} above is true 
  if $V$ is the variety generated by a finite lattice. In fact,
  for algebras generating congruence distributive varieties, a stronger
  result holds. By J{\'o}nsson's Lemma \cite[Corollary~IV.6.10]{BS:ACIU},
  every finitely generated congruence distributive variety $V$ contains only
  finitely many subdirectly irreducible algebras. Therefore $V$
  has only finitely many subvarieties and all of them are finitely generated.
  For varieties of groups, a similar result holds. Here, only
  certain subdirectly irreducible groups, namely \emph{critical} groups, are used.
  A finite group $\ab{G}$ is \emph{critical} if it does not lie in the
  variety generated by the groups in $\VarH \VarS (\ab{G})$ of size smaller
  than $\ab{G}$.
  When Oates and Powell proved that for every finite group $\ab{G}$ the 
  variety $\Var(\ab{G})$ is finitely based (i.e., determined by a finite
  set of identities, \cite{OP:IRIF}), they showed that $\Var(\ab{G})$ contains 
  only finitely many critical groups.
  Thus $\Var (\ab{G})$ has only finitely many subvarieties, and all
  of them are finitely generated.
  
In this note, we consider
 a vast generalization of groups and lattices, namely 
 those algebras that have an \emph{edge operation}
among their term functions. These edge operations were introduced
by Berman, Idziak, Markovi\'c, McKenzie, Valeriote, and Willard
in \cite{BI:VWFS} as a common generalization of Mal'cev operations
and near-unanimity operations.
A ternary operation $d$ on a set $A$ is called
a \emph{Mal'cev operation} if $d(x,x,y) = d(y,x,x) = y$ for all $x,y \in A$.
For a group, $d(x,y,z) := x  y^{-1}  z$ is a Mal'cev operation. 
 On a quasigroup $(A,\cdot,\backslash,/)$ we have the Mal'cev operation
 $d(x,y,z) := (x/(y\backslash y))\cdot(y\backslash z)$.
For $k \ge 3$, a $k$-ary operation $m$ is a
 \emph{near-unanimity operation} 
 if $m(y,x,x, \ldots, x, x) = m(x,y,x, \ldots,x, x) = \cdots = m (x,x,x, \ldots, x, y) = x$
 for all $x,y\in A$.
 On a lattice, $m (x_1, x_2, \ldots, x_k) := (x_1 \meet x_2) \join (x_1 \meet x_3) \join (x_2 \meet x_3)$
  is an example of a $k$-ary near-unanimity operation.
 For $k\geq 2$, a function $t: A^{k+1}\rightarrow A$ is a
 $k$-\emph{edge operation} if for all $x,y\in A$ we have
\[ t(y,y,x,\ldots,x) = t(y,x,y,x,\ldots,x) = x \]
 and for all $i \in \{4,\ldots, k+1\}$ and for all $x,y \in A$, we have
\[ t(x, \ldots,x, y, x, \ldots, x) = x, \text{ with } y \text{ in position } i. \]
 Hence a ternary operation $d$ is a Mal'cev operation if and
 only if $t(x,y,z) := d(y,x,z)$ is a $2$-edge operation,
 and for $k \ge 3$, a $k$-ary operation $m$ is a 
 near-unanimity operation if and only if 
 $t(x_1,\ldots,x_{k+1}) := m(x_2,\ldots,x_{k+1})$ is a $k$-edge operation.
 We say that an algebra has a Mal'cev term (near-unanimity term, edge term) if
 it has a Mal'cev operation (near-unanimity operation, edge operation) among its term functions.
 Hence every algebra with a Mal'cev term and every algebra with near-unanimity
 term has an edge term.
  The class of algebras with an edge term therefore
 contains all groups and their expansions (such as rings, vector spaces, Lie algebras \ldots),
all quasigroups, loops,
 as well as all lattices and their expansions. 
 In \cite{BI:VWFS}, we find the combinatorial characterization that a finite algebra $\ab{A}$
 has an edge term if and only if there is a polynomial $p$ with real coefficients such that for all $n \in \N$,
 $\ab{A}^n$ has at most $2^{p(n)}$ subalgebras.

 In the present paper,
 we obtain the following result on subvarieties of
 a finitely generated variety with edge term.
 \begin{thm} \label{thm:accsv1}
    Let $\ab{A}$ be a finite   algebra with an edge term.
    Then there is no infinite ascending chain $V_1 \subset V_2 \subset \cdots$
            of subvarieties of $\Var (\ab{A})$.
\end{thm}
 From this result, one can infer that property~\eqref{eq:prop} is true
 for finitely generated varieties with an edge term:
\begin{thm} \label{thm:accsv2}
    Let $\ab{A}$ be a finite algebra with an edge term.
    Then for every variety $W \subseteq \Var (\ab{A})$, there is
              a finite algebra $\ab{B} \in W$ such that
              $W = \Var (\ab{B})$.
\end{thm} 
 We will prove these theorems in Section~\ref{sec:acc}.
 Unlike for groups and lattices, there are algebras with an edge term
 that generate varieties with infinitely many subvarieties. Thus a proof 
 requires an approach that is different from the one for the classical cases.
 However, since a finitely generated variety contains, up to isomorphism, at most
  countably many finite algebras, we obtain, 
  as a consequence of Theorem~\ref{thm:accsv2}:
\begin{cor} \label{cor:number}
 Every finitely generated variety with an edge term has at most countably many subvarieties.
\end{cor}

 To illustrate that an arbitrary variety with an edge term may have uncountably many subvarieties,
 we recall two fundamental classical results:
 Evans and Neumann showed that the number of varieties of loops is continuum~\cite{EN:VGL},
 and Ol'{\v{s}}anski{\u\i} proved that there are even continuum many distinct varieties of
 groups~\cite{Ol:FBP}.

 As an application of our results to the theory of loops we obtain the following.

\begin{cor} \label{cor:loop}
 Let $\ab{A}$ be a finite loop. Then $\Var(\ab{A})$ has at most countably many subvarieties, and
 all of them are finitely generated. 
\end{cor}

  As mentioned above,
  for groups Theorem~\ref{thm:accsv2} was obtained in \cite{OP:IRIF} in the course of 
  proving that
  for every finite group $\ab{G}$, the variety $\Var (\ab{G})$ is 
  finitely based.
  In general,
  a variety $V$ of algebras of some type $\mathcal{F}$ 
   is called \emph{finitely based} if there is a finite
  set $\Phi$ of identities such that an algebra $\ab{A}$ of type
$\mathcal{F}$ lies in $V$  if and only if
$\ab{A} \models \Phi$. 
   An algebra is called 
  finitely based if the variety it generates is finitely based.
   Theorem~\ref{thm:accsv1} asserts that
  for a finite algebra $\ab{A}$ with an edge term, the class of subvarieties
  of $\Var (\ab{A})$
  has no infinite ascending chains. However, infinite descending chains
  may exist. As a corollary of the main result, we obtain
  several conditions that are equivalent to the fact that there
  are no infinite descending chains of subvarieties of $\Var (\ab{A})$.
  We say that a variety $W$ contained in a variety $V$ is \emph{finitely based 
  relative to $V$} if there is a finite set of identities $\Phi$ such that
  $W = \{ \ab{A} \in V \setsuchthat \ab{A} \models \Phi \}$. An algebra
  $\ab{B} \in V$ is finitely based relative to $V$ if $\Var(\ab{B})$ is 
  finitely based relative to $V$.
  A finite algebra $\ab{A}$ is called \emph{cardinality critical} if it
  does not lie in the variety generated by the algebras $\ab{B} \in \Var(\ab{A})$
  of size smaller than $\ab{A}$. 
  Clearly, every cardinality critical finite group is critical.
  Now combining Theorem~\ref{thm:accsv2} with known facts, we obtain
  the following corollary.
\begin{cor} \label{cor:4}
    Let $\ab{A}$ be a finite algebra with an edge term. Then
    the following are equivalent:
    \begin{enumerate}
        \item \label{it:c1} Each $\ab{B} \in \Var (\ab{A})$ is finitely based relative to
              $\Var (\ab{A})$.
        \item \label{it:c2} $\Var (\ab{A})$ has only finitely many subvarieties.
        \item \label{it:c3} $\Var (\ab{A})$ contains, up to isomorphism,
               only finitely many cardinality critical
              members.
        \item \label{it:c4} There is no infinite descending chain of subvarieties
              of $\Var (\ab{A})$.
    \end{enumerate}
\end{cor}
This corollary will be derived in Section~\ref{sec:dcc}.

There are examples of finite algebras with an edge term that satisfy none of the conditions~\eqref{it:c1}~to~\eqref{it:c4}
given in Corollary~\ref{cor:4}; one such example has been constructed in \cite{Br:TLOF}.
This example is a pointed group.
A \emph{pointed group} is an algebra $\ab{A}$ in the language 
$\mathcal{P} = \{ \cdot, \mbox{}^{-1}, 1, c \}$
such that its reduct $\ab{G} := \algop{A}{\cdot^{\ab{A}}, {\mbox{}^{-1}}^{\ab{A}}, 1^{\ab{A}}}$ is a group,
and $c$ is a nullary operation; we also write $(\ab{G}, c^{\ab{A}})$ for $\ab{A}$.
Now \cite{Br:TLOF} provides  a finite group $\ab{P}$ and an element $p \in P$ such that
the pointed group $(\ab{P}, p)$ is not finitely based
(in the language $\mathcal{P}$).
The variety generated by $\ab{Q} := \prod_{g \in P} (\ab{P}, g)$ is finitely 
based. In fact, let $\Phi$ be the finite set of identities defining the
variety generated by the group $\ab{P}$; such a set exists by the Oates-Powell-Theorem.
Then $\Phi$ also axiomatizes the variety of pointed groups generated by
$\ab{Q}$. Hence $(\ab{P}, p)$ is not finitely based relative to $\Var(\ab{Q})$,
and therefore $\ab{Q}$ does not satisfy any of the conditions~\eqref{it:c1}~to~\eqref{it:c4}
given in Corollary~\ref{cor:4}. Altogether, we have that  the class of subvarieties of
$\Var(\ab{Q})$ satisfies the ascending chain condition, but not the descending chain condition.

The main technique for proving Theorem~\ref{thm:accsv1} will be an application of the arguments that were
used in~\cite{Ai:CMCO,AMM:OTNO} to establish that
every clone on a finite set that contains an edge operation
is finitely related.
In our setting, the role of clones is taken by a suitable encoding
of the equational theory of a variety into a structure that
we will call a \emph{clonoid}. While a clonoid now represents
the equational theory of a variety, the role of the invariant relations appearing
in clone theory is taken by the algebras that lie in a variety.

\section{Preliminaries from universal algebra}

In this section, we review the relation between ``finitely generated'' and 
``no ascending chains of subvarieties'', and between
``finitely based'' and ``no descending chains of subvarieties''.

 We will first state a well-known lemma.    
For an algebra $\ab{A}$, we say that $\ab{A}$ is $k$-generated
 if it is generated by a subset $S \subseteq A$ with $|S| \le k$.
\begin{lem} \label{lem:smallalg} \label{lem:smallid}
   Let $\F$ be a type of algebras,
   let $V$ be a variety of algebras of type $\F$,  
   let $k \in \N$,
    let $\varphi := (s \approx t)$ be an identity over $\F$ that
   uses at most $k$ variables, and let  $\ab{A}$ be a $k$-generated  algebra of type $\F$.
    Then we have:
    \begin{enumerate}
          \item \label{it:1}  $V \models \varphi$ if and only if every
                $k$-generated algebra in $V$ satisfies $\varphi$.
          \item \label{it:2}  $\ab{A} \in V$ if and only if $\ab{A}$ satisfies every identity
                                with at most $k$ variables that holds in $V$.
    \end{enumerate}
\end{lem}
  For groups, item~\eqref{it:2} is explicitly given in \cite[Lemma 16.1]{Ne:VOG}.
  A variety is called \emph{locally finite} if all of its finitely generated
  members are finite. Every finitely generated variety is locally finite.
 From the lattice of subvarieties of a locally finite variety, 
  finitely generated varieties can
 be recognized in the following way. 
\begin{lem}[cf. {\cite[p.\ 370]{OV:VTMO}}] \label{lem:acc}
   Let $V$ be a locally finite variety of arbitrary type $\F$, and
   let $W$ be a subvariety of $V$. Then the following are equivalent:
   \begin{enumerate} 
      \item \label{it:a1}
                   There exists no infinite strictly
    ascending
    chain of varieties $V_1 \subset V_2 \subset V_3 \subset \cdots$ 
    with  $W := \Var (\bigcup_{i \in \N} V_i)$.
      \item \label{it:a2}
       $W$ is finitely generated.
    \end{enumerate}
\end{lem}
 Similarly, finitely based varieties will be singled out in the following
lemma.  We recall that for a set of identities $\Phi$ of type $\F$, 
  $\Mod (\Phi)$ denotes the class of all algebras of type $\F$  satisfying 
  $\Phi$.
\begin{lem}[cf. {\cite[p.\ 370]{OV:VTMO}}] \label{lem:fb}
    Let $V$ be a locally finite variety of arbitrary 
    type $\F$, and let $W$ be a subvariety of $V$.
    Then the following are equivalent:
    \begin{enumerate}
         \item 
        \label{it:l1}
            There exists no infinite strictly
    descending
    chain of varieties $V \supseteq V_1 \supset V_2 \supset V_3 \supset \cdots$ 
    with  $W := \bigcap_{i \in \N} V_i$.
         \item
        \label{it:l2} 
            There is a $k \in \N$ and a  set $\Sigma$ of identities such that
             each $\sigma \in \Sigma$ contains at most $k$ variables and
            $W = V \cap \Mod (\Sigma)$.
         \item
             \label{it:l3}
             There is a finite set $\Phi$ of identities such that
              $W = V \cap \Mod (\Phi)$.
     \end{enumerate}
\end{lem}

   Finally, we need a lemma on cardinality critical algebras.
 We recall that  a finite algebra $\ab{B}$ is called
     \emph{cardinality critical} if 
     $\ab{B} \not\in \Var ( \{ \ab{C} \setsuchthat \ab{C} \in \Var(\ab{B}),
                                                |\ab{C}| < | \ab{B} | \})$.
    We have:
    \begin{lem} \label{lem:cc}
        Let $V$ be a locally finite variety. Then each subvariety
        $W$ of $V$ is generated by the cardinality critical algebras
        it contains.
    \end{lem}
     \emph{Proof:}
      Let $W_1 := \Var (\{ \ab{B} \in W \setsuchthat \ab{B} \text{ is finite and cardinality critical}\})$.
    Since $W$ and  $W_1$ are both locally finite, it suffices to show that they have the
   same finite members. Seeking a contradiction, we let $\ab{C} \in W$ be finite and of minimal size
    with $\ab{C} \not\in W_1$. Then $\ab{C}$ is not cardinality critical, hence
     it lies in the variety generated by $X := \{ \ab{D} \setsuchthat |\ab{D}| < |\ab{C}|,
                                            \ab{D} \in \Var (\ab{C}) \}$. By the minimality
     of $\ab{C}$, we have $X \subseteq W_1$, and thus $\ab{C} \in W_1$, a contradiction. \qed

\section{Preliminaries from order theory}

The techniques that we use from order theory have been used
in a similar way in \cite{Ai:CMCO}, and in an almost
identical way as they are used in the present note
in \cite{AMM:OTNO}. 
To keep our presentation self-contained, we report these results
with only a slight modification in notation.

A partially ordered set $\algop{X}{\le}$ satisfies
the \emph{descending chain condition} $\DCC$ if
there are no
infinite descending chains $x_1 > x_2 > \cdots$ in $X$,
and it satisfies \emph{the ascending chain condition} $\ACC$ if
there is no infinite ascending chain $x_1 < x_2 < \cdots$ in $X$.
A subset $Y$ of $X$ is an \emph{antichain} if for
all $y_1, y_2 \in Y$ with $y_1 \le y_2$, we have
$y_1 = y_2$. 
The partially ordered set $\algop{X}{\le}$ is
\emph{well partially ordered} if it satisfies the $\DCC$
and has no infinite antichains.
For a partially ordered set $\algop{X}{\le}$, a subset $Y$ of $X$
is \emph{upward closed} if for all $y \in Y$ and $x \in X$ with
$y \le x$, we have $x \in Y$.
The set of upward closed subsets of $X$ is denoted by
$\U ({X}, {\le})$. If $\algop{X}{\le}$ is well partially
ordered, then $\algop{\U ({X},{\le})}{\subseteq}$
satisfies the $\ACC$ (cf. \cite[Theorem~1.2]{Mi:BWAB2}). 

For $A = \{1, 2, \ldots, t \}$,  
we will use the lexicographic ordering
on $A^n$. 
For $\vb{a} = (a_1,\ldots, a_n)$ and
$\vb{b} = (b_1, \ldots, b_n)$, we say
$\vb{a} \lelex \vb{b}$ if
\begin{multline*}
     (\exists i \in \{1,\ldots,n\} :
       a_1 = b_1 \wedge \ldots \wedge
       a_{i-1} = b_{i-1} \wedge a_i < b_i) \text{ or }
      \\ (a_1, \ldots,a_n) = (b_1,\ldots, b_n).
\end{multline*}

For every finite set $A$,
we let $A^+$ be the set
$\bigcup \{A^n \setsuchthat n \in \N\}$. (We use
 $\N$ for the set of natural numbers $\{1,2,3, \ldots \}$.)
We will now introduce an order relation on $A^+$.
For $\vb{a} = (a_1,\ldots, a_n) \in A^+$ and $b \in A$,
we define the \emph{index of the first occurrence of $b$ in $\vb{a}$},
$\FO (\vb{a}, b)$, by
$\FO (\vb{a}, b) := 0$ if $b \not\in \{a_1,\ldots, a_n \}$,
and $\FO (\vb{a}, b) := \min \{ i \in \{1,\ldots, n\} \setsuchthat
a_i = b\}$ otherwise.
\begin{de}[{\cite[Definition~3.1]{AMM:OTNO}}]
    Let $\vb{a} = (a_1,\ldots, a_m)$ and $\vb{b} = (b_1,\ldots, b_n)$
    be in $A^+$. We say
    $\vb{a} \lea \vb{b}$ if there is an injective and
    increasing function $h :\{1,\ldots, m\} \to \{1,\ldots,n\}$
    such that
    \begin{enumerate}
         \item for all $i \in \{1,\ldots, m\}$ : 
                       $a_i = b_{h(i)}$,
         \item $\{a_1,\ldots, a_m\} = \{b_1,\ldots, b_n\}$,
         \item for all $c \in \{a_1,\ldots, a_m\}$:
                $h (\FO (\vb{a}, c)) = \FO (\vb{b}, c)$.
    \end{enumerate}
We will call such an $h$ a function \emph{witnessing $\vb{a} \lea \vb{b}$}.
\end{de}
Informally, $\vb{a} \lea \vb{b}$ if and only if  $\vb{a}$ and $\vb{b}$ contain
the same set of letters, and $\vb{b}$ can be obtained from $\vb{a}$ by 
certain insertions, where
the insertion $\ab{x}\ab{y} \to \ab{x} l \ab{y}$ is allowed only
           if the letter $l$ appears in the word $\ab{x}$.

\begin{pro}[{cf.\ \cite[Lemma~3.2]{AMM:OTNO}}]\label{pro:leawpo}
  Let $A$ be a finite set.
  Then $\algop{A^+}{\lea}$ is well partially ordered,
  and $\algop{\U({A^+},{\lea})}{\subseteq}$ satisfies the $\ACC$.
\end{pro}
\emph{Proof:}
    In \cite[Lemma~3.2]{AMM:OTNO}, it is proved that
    $\algop{A^+}{\lea}$ is well partially ordered.
    Since $\algop{A^+}{\lea}$ is well partially ordered,
    \cite[Theorem~1.2]{Mi:BWAB2} yields that
    $\algop{\U({A^+}, {\lea})}{\subseteq}$ satisfies the $\ACC$. \qed

We will now give a slight modification of \cite[Definition~3.3]{AMM:OTNO}
and \cite[Lemma~3.4]{AMM:OTNO}.
\begin{de}[{cf. \cite[Definition~3.3]{AMM:OTNO}}]  Let $A$ be a finite set, 
    let $\vb{a} = (a_1,\ldots,a_m) \in A^m$, 
    $\vb{b} = (b_1, \ldots, b_n) \in A^n$
    be such that $\vb{a} \lea \vb{b}$, and let $h$ be a function
    from $\{1,\ldots, m\} \to \{1,\ldots, n\}$ witnessing
    $\vb{a} \lea \vb{b}$.
    We will now define a function $\Tab : \{1, \ldots, n\} \to \{1, \ldots, m\}$.
    If $j \in \ran (h)$, then
    $\Tab (j)$ is defined by
    \[
       \Tab (j) := h^{-1} (j).
    \]
      If $j \not\in \ran (h)$,
   then
   \[
      \Tab (j) := i,
   \]
    where $i$ is minimal in $\{1,\ldots, m\}$ with
   $a_i = b_j$.
  \end{de}

 \begin{lem}[{\cite[Lemma~3.4]{AMM:OTNO}}] \label{lem:co}
   Let $t \in \N$, let $A = \{ 1, 2, \ldots, t \}$, let
   $\ab{a} \in A^m$, $\ab{b} \in A^n$, and let
   $h : \{1, \ldots, m\} \to \{1,\ldots, n\}$ be a function
   witnessing $\ab{a} \lea \ab{b}$.
   Let $\vb{c} \in A^m$ be such that
   $\vb{c} \ltlex \vb{a}$.
   Then we have
   \begin{enumerate} 
      \item \label{it:c1} $\langle a_{\Tab (j)} \setsuchthat j \in \{1,\ldots, n\} \rangle = \vb{b}$,
      \item \label{it:c2} $\langle c_{\Tab (j)} \setsuchthat j \in \{1,\ldots, n\} \rangle \ltlex \vb{b}$.
   \end{enumerate}
\end{lem}

 The following proof is similar to the proof of~\cite[Lemma~3.4]{AMM:OTNO} with some notational
 changes.

\emph{Proof:}
    For proving \eqref{it:c1}, let $j \in \{1, \ldots, n \}$.
    We have to prove that $a_{\Tab (j)} = b_j$. If $j \in \ran (h)$, then
    $a_{\Tab (j)} = a_{h^{-1} (j)} = b_{h (h^{-1} (j))} = b_j$.
    If $j \not\in \ran (h)$, then $\Tab (j)$ has the property
     $a_{\Tab (j)} = b_j$. 
    
     For proving~\eqref{it:c2},
    let $k$ be the index of the first place in which
    $\vb{c}$ differs from $\vb{a}$.
    Hence $\vb{c} = (a_1,\ldots, a_{k-1}, c_k, c_{k+1},\ldots)$,
    $\vb{a} = (a_1,\ldots, a_{k-1}, a_k, a_{k+1},\ldots)$, and
    $c_k < a_k$.
    
     We first show that 
     for all $j < h(k)$, we have
    $c_{\Tab(j)}=  a_{\Tab(j)}$.
     If $j$ is in the range of $h$, we have
     $c_{\Tab (j)}  = c_{h^{-1} (j)}$ and
     $a_{\Tab (j)}  = a_{h^{-1} (j)}$.
     Since $h( h^{-1} (j)) < h(k)$, we have, by the monotonicity of $h$,
     $h^{-1} (j) < k$. Thus $c_{h^{-1} (j)} = a_{h^{-1} (j)}$, since $k$ is the first index
     at which $\vb{c}$ and $\vb{a}$ differ.
     We now consider the case that
     $j$ is not in the range of $h$. Since
     $\{ b_1,\ldots, b_n\} = \{a_1,\ldots, a_m\}$,
     there is an $i_1 \in \{1, \ldots, m\}$ such that
     $a_{i_1} = b_j$. 
     Let $i_2 := \FO (\vb{a}, a_{i_1})$.
     By the definition of
     $\Tab$, we have $\Tab (j) := i_2$.
     Now let $j_2 := h(i_2)$.
     Then we have $j_2 = \FO (\vb{b}, b_j)$, and therefore
     $j_2 \le j$. Hence $j_2 < h(k)$, and thus
     $i_2 < k$. Therefore $c_{i_2} = a_{i_2}$, and thus
     $c_{\Tab (j)} = a_{\Tab (j)}$.
     
     Since $\Tab (h(k)) = k$, we have
     $a_{\Tab (h(k))} = a_k$ and $c_{\Tab (h(k))} = c_k$, and therefore
     $\langle c_{\Tab (j)} \setsuchthat j \in \{1,\ldots, n\} \rangle \ltlex
     \langle a_{\Tab (j)} \setsuchthat j \in \{1,\ldots, n\} \rangle = \vb{b}$.\qed

\section{Clonoids}
 In this section, we will introduce \emph{clonoids}. They are sets
 of finitary functions from a set $A$ into the universe of an algebra $\ab{B}$ that
 are closed under the operations of $\ab{B}$ and under manipulation of arguments.

\begin{de}  \label{de:clonoid}
  Let $\ab{B}$ be an algebra, and let $A$ be a nonempty set.
   For a subset $\cb{C}$ of $\bigcup_{n \in \N} B^{A^n}$ and $k \in \N$, we let
   $\cb{C}^{[k]} := \cb{C} \cap B^{A^k}$.
   We call $\cb{C}$ a \emph{clonoid with source set $A$ and target algebra $\ab{B}$} if 
   \begin{enumerate}
       \item for all $k \in \N$: $\cb{C}^{[k]}$ is a subuniverse
             of $\ab{B}^{A^k}$, and
       \item \label{it:dc1} 
             for all $k,n  \in \N$, for all
             $(i_1,\ldots, i_k) \in \{1,\ldots, n\}^k$, and for
             all $c \in \cb{C}^{[k]}$, the function
             $c' : A^n \to B$ with $c' (a_1,\ldots, a_n) := 
             c (a_{i_1}, \ldots, a_{i_k})$ lies in $\cb{C}^{[n]}$.
   \end{enumerate}
\end{de}
We note that every clone $\cb{C}$ on a set $A$ is
a clonoid with source set $A$ and target algebra $\algop{A}{\cb{C}}$
 (see~\cite[Definition 4.2]{MMT:ALVV}).

Let $\ab{A}$ be an algebra, let
$m \in \N$, and let $F$ be a subuniverse
of $\ab{A}^m$. For $i \in \{1,\ldots, m\}$, we define
the relation $\varphi_i (F)$ on $A$ by
\begin{multline}
    \varphi_i (F) := \{ (a_i, b_i) \setsuchthat \\
                        (a_1,\ldots, a_m) \in F,
                        (b_1,\ldots, b_m) \in F,
                        (a_1, \ldots, a_{i-1}) =
                        (b_1, \ldots, b_{i-1}) \}.
\end{multline}
An element of $\varphi_i (F)$ is also called a \emph{fork} of $F$ at place $i$.
We will now encode algebras by their forks. As a matter of fact,
the sequence of forks does not contain the complete information about a subuniverse
of $\ab{A}^m$,
but if $\ab{A}$ has an edge term, then 
the information is helpful in distinguishing between two subalgebras
of $\ab{A}^m$
such that one is contained in the other.
For a subset $F$ of $A^m$ and $T \subseteq \{1,\ldots, m\}$, we will use
$\pi_{T} (F)$ for the projection to the components with indices in $T$. Formally,
seeing an element $\ab{a}$ as a function from $\{1,\ldots, m\}$ to $A$, we have
$\pi_{T} (F) := \{ \ab{a}|_T \setsuchthat \ab{a} \in F \}$.
\begin{lem}[{\cite[Lemma~4.1]{AMM:OTNO}}] \label{lem:FG}
   Let $k,m \in \N, k\geq 2$, and let $\ab{A}$ be an algebra with a $k$-edge term. 
   Let $F,G$ be subuniverses of $\ab{A}^m$ with $F \subseteq G$.
   Assume that $\varphi_i (G) = \varphi_i (F)$ for all $i \in \{1,\ldots, m\}$ and
   $\pi_T(F) = \pi_T(G)$ for all $T\subseteq\{1,\dots,m\}$ with $|T|<k$.  
   Then $F=G$.
\end{lem}

For each $k \in \N$, we will now encode the forks of the algebra of $k$-ary
functions in a clonoid  in a way that runs completely
parallel to the encoding used for clones in \cite{Ai:CMCO} and~\cite{AMM:OTNO}.

Let $A$ be the set $\{1, \ldots, t\}$, let $\ab{B}$ be an algebra,
and let $\cb{C}$ be a clonoid with source set $A$ and target algebra $\ab{B}$.
For $n \in \N$ and $\vb{a} \in A^n$, we define a binary 
relation $\PhiX (\cb{C}, \vb{a})$ on $B$ by
\[
    \PhiX (\cb{C}, \vb{a}) :=
       \{ (f(\vb{a}), g(\vb{a})) \setsuchthat
          f,g \in \cb{C}^{[n]},
          \forall \vb{c} \in A^n: 
            \vb{c} \ltlex \vb{a} 
             \Rightarrow f (\vb{c}) = g(\vb{c}) \}.
\]
 Hence the elements of $\PhiX (\cb{C}, \vb{a})$ are the forks of $\cb{C}^{[n]}$ at $\vb{a}$.
We have the following lemma.
\begin{lem}[{cf. \cite[Lemma~5.2]{AMM:OTNO} and
             \cite[Lemma~4.1]{Ai:CMCO}}] \label{lem:embeddingPhi}
    Let $m,n, t \in \N$, let $\ab{B}$ be an algebra, and let $\cb{C}$ be a 
    clonoid with source set $A =\{1,2,\ldots, t\}$ and target algebra
    $\ab{B}$.
    Let $\vb{a} \in A^m$, $\vb{b} \in A^n$ such that
    $\vb{a} \lea \vb{b}$.
    Then $\PhiX (\cb{C}, \vb{b}) \subseteq \PhiX (\cb{C}, \vb{a})$.
\end{lem}
\emph{Proof:}
   Let $(x,y) \in \PhiX (\cb{C}, \vb{b})$. Then there are
   $f,g \in \cb{C}^{[n]}$ such that $x = f(\vb{b})$,
   $y = g(\vb{b})$, and $f(\vb{c}) = g (\vb{c})$ for all
   $\vb{c} \in A^n$ with $\vb{c} \ltlex \vb{b}$.
   Let $h$ be a function from $\{1,\ldots, m\}$ to $\{1,\ldots,n\}$
   witnessing $\vb{a} \lea \vb{b}$.
   Now we define functions $f_1$ and $g_1$ from $A^m$ to $A$ by
 \[
     \begin{array}{rcl}
        f_1 (\vb{x}) & := & f (\langle x_{\Tab(j)} \setsuchthat j \in \{1, \ldots, n \} \rangle) \\
        g_1 (\vb{x}) & := & g (\langle x_{\Tab(j)} \setsuchthat j \in \{1, \ldots, n \} \rangle) \\
     \end{array}
 \]
  for $\vb{x}\in A^m$. By the closure properties of a clonoid, $f_1$ and $g_1$ lie in the clonoid $\cb{C}$.
 
We will now show that $(f_1 (\vb{a}), g_1 (\vb{a}))$ is an element of 
$\PhiX (\cb{C}, \vb{a})$. To this end,
 let $\vb{c} \in A^m$ be such that $\vb{c} \ltlex \vb{a}$.
 Then we have
 \(
     f_1 (\vb{c}) 
     =
           f (\langle c_{\Tab(j)} \setsuchthat j \in \{1, \ldots, n \} \rangle).
 \)
 Since $\vb{c} \ltlex \vb{a}$, Lemma~\ref{lem:co} yields
$\langle c_{\Tab(j)} \setsuchthat j \in \{1, \ldots, n \} \rangle 
 \ltlex \vb{b}$.
 Hence we have
 $f_1(\vb{c}) = f (\langle c_{\Tab(j)} \setsuchthat j \in \{1, \ldots, n \} \rangle) = 
  g (\langle c_{\Tab(j)} \setsuchthat j \in \{1, \ldots, n \} \rangle  = g_1 (\vb{c})$.
From this, we obtain $(f_1 (\vb{a}), g_1 (\vb{a})) \in \PhiX (\cb{C}, \vb{a})$.
Since by Lemma~\ref{lem:co}~\eqref{it:c1},
      $(f_1 (\vb{a}), g_1 (\vb{a})) = (f (\vb{b}), g (\vb{b}))
        = (x,y)$, we obtain $(x,y) \in \PhiX (\cb{C}, \vb{a})$.
\qed
                                  
\begin{de}[{cf. \cite[Definition~4.2]{Ai:CMCO}}]
   Let $t \in \N$, let $\cb{C}$ be a clonoid with source set 
  $A = \{1,2, \ldots, t \}$
   and target algebra $\ab{B}$, and let
   and let $\alpha \subseteq B \times B$.
   We define a subset $\Psi (\cb{C}, \alpha)$ of $A^+$ by
   \(
      \Psi (\cb{C}, \alpha) :=
        \{ \vb{a} \in A^+ \setsuchthat
           \PhiX (\cb{C}, \vb{a}) \subseteq \alpha \}.
   \)
\end{de}
 Hence $\Psi (\cb{C}, \alpha)$ is the set of those places at which all forks 
 from $\cb{C}$ lie in $\alpha$.

\begin{lem}[{cf. \cite[Lemma~4.3]{Ai:CMCO}}] \label{lem:upwardclosed}
    Let $t \in \N$, let $\cb{C}$ be a 
    clonoid with source set $A =\{1,2,\ldots, t\}$ and target
    algebra $\ab{B}$,
     and let $\alpha \subseteq B \times B$.
    Then $\Psi (\cb{C}, \alpha)$ is an upward closed subset
   of $\algop{A^+ }{\lea}$.
\end{lem}
\emph{Proof:}
    Let $\vb{a} \in \Psi(\cb{C}, \alpha)$, and let
    $\vb{b} \in A^+$ such that $\vb{a} \lea \vb{b}$.
    Since $\vb{a} \in \Psi (\cb{C}, \alpha)$, we have
    $\PhiX (\cb{C}, \vb{a}) \subseteq \alpha$.
    By Lemma~\ref{lem:embeddingPhi}, we have
    $\PhiX (\cb{C}, \vb{b}) \subseteq \PhiX (\cb{C}, \vb{a})$.
     Therefore, $\PhiX (\cb{C}, \vb{b}) \subseteq \alpha$,
    and thus $\vb{b} \in \Psi (\cb{C}, \alpha)$. \qed

\section{Chains of Clonoids}

 The next lemma allows us to decide when two comparable clonoids are equal.

\begin{lem} \label{lem:CD}  Let $k, t \in \N$ with $k\geq 2$, let $A := \{1,2,\ldots, t\}$,
            let $\ab{B}$ be an algebra with $k$-edge term,
            let $(\mathcal{C}, \subseteq)$ be a linearly ordered set
            of clonoids with source set $A$ and target algebra $\ab{B}$, and
            let $\cb{C}, \cb{D} \in \mathcal{C}$. Then the following are
            equivalent:
            \begin{enumerate}
                \item \label{it:l1}
                      $\cb{C} \subseteq \cb{D}$.
                \item \label{it:l2}
                      $\cb{C}^{[t^{k-1}]} \subseteq \cb{D}^{[t^{k-1}]}$ and
                      for all $\alpha \subseteq B \times B$, we
                      have $\Psi (\cb{D}, \alpha) \subseteq \Psi (\cb{C}, \alpha)$.     
            \end{enumerate}
\end{lem}

\emph{Proof:}
   \eqref{it:l1}$\Rightarrow$\eqref{it:l2}:
      Let $\alpha \subseteq B\times B$, and let $\vb{a} \in
       \Psi (\cb{D}, \alpha)$.
       Then $\PhiX (\cb{D}, \vb{a}) \subseteq \alpha$. Since
       $\cb{C} \subseteq \cb{D}$, we therefore have
       $\PhiX (\cb{C}, \vb{a}) \subseteq \alpha$, and thus
       $\vb{a} \in \Psi (\cb{C}, \alpha)$.

      \eqref{it:l2}$\Rightarrow$\eqref{it:l1}:
          Since $\mathcal{C}$ is a linearly ordered set of clonoids,
         we either have $\cb{C} \subseteq \cb{D}$ or
         $\cb{D} \subseteq \cb{C}$. In the first case, there is nothing
         to prove, so we assume $\cb{D} \subseteq \cb{C}$.

        We will now prove that in this case, we have
        $\cb{C}^{[n]} \subseteq \cb{D}^{[n]}$ for all $n \in \N$.
        Let us first assume  $n \le t^{k-1}$. In this case, the result is 
        a consequence of the closure property of clonoids that
        is given in item~\eqref{it:dc1}
        of Definition~\ref{de:clonoid}: 
 indeed, if $f \in \cb{C}^{[n]}$, then
        $f' : (x_1,\ldots, x_{t^{k-1}}) \mapsto f(x_1,\ldots, x_n)$ is an
        element of $\cb{C}^{[t^{k-1}]}$, and hence of $\cb{D}^{[t^{k-1}]}$.
        Hence $g : (x_1,\ldots, x_n) \mapsto f'(x_1,\ldots, x_n, x_n, \ldots, x_n)$ lies
        in $\cb{D}^{[n]}$, and since $f=g$, we obtain $f \in \cb{D}^{[n]}$.  
         Let us now assume $n > t^{k-1}$. We consider $F := \cb{C}^{[n]}$ and
 $G := \cb{D}^{[n]}$ as subuniverses of $\ab{B}^{A^n}$, and we will employ
        Lemma~\ref{lem:FG} to show $F = G$. To this end, we
        first show that for all $T\subseteq A^n$ with $|T|< k$ we have
        \begin{equation} \label{eq:pcd}
           \pi_T (F) \subseteq \pi_T (G).
        \end{equation}
 Let $f\in F$, and let $T\subseteq A^n$ with $|T|< k$.

 For $\vb{a} = (a_1,\dots,a_n)$ in $T$, let 
 $\mu_{\vb{a}} := \{(i,j)\in\{1,\dots,n\}^2 \setsuchthat a_i = a_j \}$.
 Then $\mu_{\vb{a}}$ is a partition of $\{1,\dots,n\}$ into at most $t$ blocks.
 Let $\mu := \bigcap_{\vb{a}\in T} \mu_{\vb{a}}$, and let $q$ be the number of equivalence
 classes modulo $\mu$. Clearly $q \leq t^{k-1}$.
 Put differently, let $M$ be the $(|T| \times n)$-matrix with entries from $A$ whose rows are the
 elements of $T$, and let $(i,j) \in \mu$ if the $i$-th column of $M$ is equal to the
 $j$-th column. Since there are at most $|A|^{|T|}$ different column vectors,
 $\mu$ has at most $|A|^{|T|} \le t^{k-1}$ equivalence classes.

 Let $S_1,\dots,S_q$ be the
 equivalence classes modulo $\mu$, and let
 $r_1 \in S_1, \dots, r_q \in S_q$ be representatives of these classes. 
 For $j\in \{1,\dots,n\}$, let $c(j) \in \{1,\dots,q\}$ be the unique element such that $j\in S_{c (j)}$.
 By the definition of a clonoid,
 $f_\mu: A^q\to B, (x_1,\dots,x_q) \mapsto f(x_{c (1)},\dots,x_{c (n)})$ is in $\vb{C}^{[q]}$.
 Since $q \leq t^{k-1}$, it follows that $f_\mu\in\vb{D}^{[q]}$ as well.
 Hence the function $g\colon A^n\to B, (x_1,\dots,x_n) \mapsto f_\mu(x_{r_1},\dots,x_{r_q})$
 is in $\vb{D}^{[n]}$.
 Let $\vb{a} = (a_1,\dots,a_n)\in T$. Then
 $g(a_1,\dots,a_n) =  f_\mu(a_{r_1},\dots,a_{r_q}) = f (a_{r_{c(1)}}, \ldots, a_{r_{c(n)}})$.
 Now for each $j \in \{1,\ldots, n\}$, we have
 $r_{c(j)} \in S_{c(j)}$ and $j \in S_{c(j)}$, thus
 $(r_{c(j)}, j) \in \mu$ and therefore $a_{r_{c(j)}} = a_j$.
 Therefore, $f (a_{r_{c(1)}}, \ldots, a_{r_{c(n)}}) = f(a_1,\ldots, a_n)$.  
 Hence $\pi_T(f) = \pi_T(g)$, and thus $\pi_T (f) \in \pi_T (G)$, which proves~\eqref{eq:pcd}.

 Next we show that for all $\vb{a} \in A^n$,
        \begin{equation} \label{eq:cd}
           \PhiX (\cb{C}, \vb{a}) \subseteq \PhiX (\cb{D}, \vb{a}).
        \end{equation}
         In order to prove~\eqref{eq:cd}, we fix $\vb{a} \in A^n$.
         We obviously have $\vb{a} \in \Psi (\cb{D}, \PhiX(\cb{D}, \vb{a}))$.
         By the assumption~\eqref{it:l2}, we therefore have
         $\vb{a} \in \Psi (\cb{C}, \PhiX(\cb{D}, \vb{a}))$
 and~\eqref{eq:cd} follows from the definition of $\Psi$.

 Now from~\eqref{eq:pcd},
          \eqref{eq:cd}, and Lemma~\ref{lem:FG}
          (with $m := |A|^n$), we obtain $\cb{C}^{[n]} =
         \cb{D}^{[n]}$. \qed

        Let $2^{B \times B}$ denote
        the power set of $B \times B$.
        On the set $\ob{U} := \U( {A^+},{\lea})^{2^{B \times B}}$,
        we define an order as follows:
        for $\ab{S} = 
        \langle S (\alpha) \setsuchthat \alpha \subseteq B \times B \rangle$ 
      and $\ab{T} =
         \langle T (\alpha) \setsuchthat \alpha \subseteq B \times B \rangle
         \in \ob{U}$, we say $\ab{S} \pmb{\le} \ab{T}$ 
        if $S(\alpha) \subseteq  T (\alpha)$ for all $\alpha \subseteq B \times B$.
         Hence $\algop{\ob{U}}{\pmb{\le}}$ is isomorphic to
        the $2^{|B|^2}$-fold direct product of $\algop{\U ({A^+}, {\lea} )}{\subseteq}$.
       Therefore, it follows from Proposition~\ref{pro:leawpo} that
       $(\ob{U}, \pmb{\le})$ satisfies the $\ACC$.          

      Now, as a corollary of Lemma~\ref{lem:CD}, we obtain:
    \begin{lem} \label{lem:R} 
                 Let $k,t \in \N$ with  $k\geq 2$, let $A := \{1,2,\ldots,t\}$, and  let $\ab{B}$ be
            an algebra with a $k$-edge term.
            Let $(\mathcal{C}, \subseteq)$ be a linearly ordered set
            of clonoids with source set $A$ and target algebra $\ab{B}$,
let $\ob{T}$ denote the power set of $B^{A^{t^{k-1}}}$, and let $\ob{U} := \U ( {A^+}, {\lea})^{2^{B \times B}}$.
            Let $R : \mathcal{C} \to \ob{T}\times\ob{U}$ be defined
            by
            \[
                R (\cb{C}) := (\cb{C}^{[t^{k-1}]},\langle \Psi (\cb{C}, \alpha) \setsuchthat
                                      \alpha \subseteq B \times B \rangle)
            \]
            for $\cb{C} \in \mathcal{C}$.
            Then $R$ is injective, and for all $\cb{C}, \cb{D} \in \mathcal{C}$ with
            $\cb{C} \subseteq \cb{D}$, we have $R(\cb{C})_1 \subseteq R(\cb{D})_1$ and
            $R (\cb{D})_2 \pmb{\le} R (\cb{C})_2$.
     \end{lem}

\begin{thm} \label{thm:dcc}
   Let $\ab{B}$ be a finite algebra with an edge term, and
   let
   $A$ be a nonempty finite set.
   Let $\mathcal{C} := \{ \cb{C} \setsuchthat
    \cb{C}$ is clonoid with source set $A$ and target algebra $\ab{B} \}$.
   Then $(\mathcal{C}, \subseteq)$ satisfies the descending chain condition.
\end{thm}

 \emph{Proof:} Let $\mathcal{C}$ be an infinite descending chain
  of clonoids with source set $A$ and target algebra $\ab{B}$.
         From these clonoids, Lemma~\ref{lem:R} produces
        an infinite ascending chain in
        $\algop{\U ({A^+}, {\lea})}{\subseteq}^{2^{B \times B}}$, which contradicts
        Proposition~\ref{pro:leawpo}. \qed

\section{Varieties} \label{sec:acc}
For an algebra $\ab{A}$, we will introduce a clonoid $\cb{C}$ with source set
$A$ and target algebra $\ab{A} \times \ab{A}$ such that the equational theory
of each subvariety of $\Var (\ab{A})$ corresponds to a ``subclonoid'' of $\cb{C}$.

\begin{de} 
   Let $\ab{A}$ be an algebra, and let $W$ be a subvariety of $\Var (\ab{A})$.
   By $\Th_{\ab{A}} (W)$ we denote the clonoid with source set $A$ and target
   algebra $\ab{A} \times \ab{A}$ defined by
      \begin{multline}
      \Th_{\ab{A}} (W) := \{ (a_1,\ldots, a_k) \mapsto \VecTwo{s^{\ab{A}} (\vb{a})}{t^{\ab{A}} (\vb{a})} \setsuchthat
                       k \in \N,  \\ 
                        s, t \text{ are $k$-variable terms in the language of } \ab{A} \text{ with }
                        W \models s \approx t \}.
   \end{multline}
\end{de}
   We note that $\Th_{\ab{A}} (W)$ clearly satisfies the closure properties from
   the definition of a clonoid.

\begin{lem} \label{lem:galois}
   Let $\ab{A}$ be an algebra, and let $W_1$ and $W_2$ be subvarieties
   of $\Var(\ab{A})$. Then we have:
     \begin{equation} \label{eq:galois}
        W_1 \subseteq W_2 \text{ if and only if } \Th_{\ab{A}} (W_2) \subseteq \Th_{\ab{A}} (W_1).
   \end{equation}
\end{lem}    
\emph{Proof:}
       For the ``only if''-direction, we observe that
      for all terms $s, t$ in the language of $\ab{A}$ with
      $W_2 \models s \approx t$, we have $W_1 \models s \approx t$.
      This implies that $\Th_{\ab{A}} (W_2) \subseteq \Th_{\ab{A}} (W_1)$.
      For the ``if''-direction, we show that every identity 
      that is satisfied in $W_2$ is also satisfied in $W_1$. To this end,
      we let $k \in \N$, and we let $s,t$ be two $k$-variable terms such that
      $W_2 \models s \approx t$. Then the $k$-ary function $\vb{a} \mapsto \VecTwo {s^{\ab{A}} (\ab{a})}
      {t^{\vb{A}} (\ab{a})}$ lies in $\Th_{\ab{A}} (W_2)$, and therefore, by the assumption,
      also in $\Th_{\ab{A}} (W_1)$. Hence there exist $k$-ary terms $s_1, t_1$ such that
      $s^{\ab{A}} = s_1^{\ab{A}}$, $t^{\ab{A}} = t_1^{\ab{A}}$, and 
      $W_1 \models s_1 \approx t_1$. Since $W_1 \subseteq \Var (\ab{A})$,
      we have $W_1 \models s \approx s_1$ and $W_1 \models t \approx t_1$.
      Hence $W_1 \models s \approx t$.
      Therefore, every identity satisfied in $W_2$ is satisfied in $W_1$, and thus
      $W_1 \subseteq W_2$. \qed

  We can now give the proof of the main results:

\emph{Proof of Theorem~\ref{thm:accsv1}:}
     Let
      $W_1 \subset W_2 \subset W_3 \subset \ldots$ be a strictly increasing sequence of subvarieties of $\Var(\ab{A})$.
      Then by Lemma~\ref{lem:galois}, $\Th_{\ab{A}} (W_1) \supset \Th_{\ab{A}} (W_2) \supset \Th_{\ab{A}} (W_3) \supset \ldots$
      is an infinite strictly decreasing subsequence of clonoids with source set $A$ and target algebra $\ab{A} \times \ab{A}$,
      contradicting Theorem~\ref{thm:dcc}. Hence  the subvarieties of $\Var(A)$ satisfy
      the ascending chain condition. \qed

\emph{Proof of Theorem~\ref{thm:accsv2}:} 
     The result follows from Theorem~\ref{thm:accsv1} and
      the implication \eqref{it:a1}$\Rightarrow$\eqref{it:a2}
       of Lemma~\ref{lem:acc}. \qed

\section{On the descending chain condition for subvarieties} \label{sec:dcc}

    In this section, we will investigate the relation of Theorem~\ref{thm:accsv2}
    to the finite basis property.
   For varieties $V$ and $W$ of the same type, we write
 $V \lcover W$ if $V \subset W$ and there
 is no variety $V'$ with $V \subset V' \subset W$, and we say that $V$ is a \emph{subcover}
  of $W$.  
    We will first prove a lemma that states that in the lattice of all varieties of a given type,
  every finitely generated variety has only finitely many subcovers.
   This property resembles the fact that every finitely related clone
   on a finite set has only finitely many covers \cite[p.\ 93, 4.1.3]{PK:FUR}.
\begin{lem} \label{lem:finitelybranching}
   Let $\ab{A}$ be a finite algebra. Then we have:
   \begin{enumerate}
     \item \label{it:co1} There are only finitely many varieties $W$ with $W \lcover \Var (\ab{A})$.
     \item \label{it:co2} For each variety $U \subset \Var (\ab{A})$ there is a variety $W$ with
           $U \subseteq W \lcover \Var (\ab{A})$.
   \end{enumerate}
\end{lem}
\emph{Proof:}
     Let $k \in \N$ be such that $\ab{A}$ is generated by $k$ elements.
     We first show that for all varieties $W \subseteq \Var (\ab{A})$ we have:
     \begin{equation} \label{eq:WVB}
        \text{If $\Th_{\ab{A}}^{[k]} (W) = \Th_{\ab{A}}^{[k]} (\Var (\ab{A}))$, then
        $W = \Var (\ab{A})$.}
     \end{equation}
     To this end, we show $\ab{A} \in W$. By Lemma~\ref{lem:smallid},
     it is sufficient to check that every identity $s \approx t$ of $W$ with at most $k$
     variables is satisfied by $\ab{A}$ .
     Letting $s^{\ab{A}}$ and $t^{\ab{A}}$ be the $k$-ary term functions induced by $s$ and $t$,
     we see that $(s^{\ab{A}}, t^{\ab{A}}) \in \Th_{\ab{A}}^{[k]} (W)$. 
     Therefore $(s^{\ab{A}}, t^{\ab{A}}) \in \Th_{\ab{A}}^{[k]} (\Var (\ab{A}))$.
     This implies that there are $k$-ary terms $s_1, t_1$ such that
     $(s_1^{\ab{A}}, t_1^{\ab{A}}) = (s^{\ab{A}}, t^{\ab{A}})$ and 
     $\Var (\ab{A}) \models s_1 \approx t_1$. From this, we get
     $s^{\ab{A}} = t^{\ab{A}}$, and hence $\ab{A} \models s \approx t$.
     This completes the proof of~\eqref{eq:WVB}.
   
  Let $R$ be a set of representatives of the $k$-ary terms in the language of $\ab{A}$ modulo
     the equivalence relation defined by $s \sim_{\ab{A}} t \Leftrightarrow
     s^{\ab{A}} = t^{\ab{A}}$.
     Let \[S := \{ \Var (\ab{A}) \cap \Mod (s \approx t) \setsuchthat
                    s,t \in R, s \not\sim_{\ab{A}} t \}. \]
Note that $R$ and hence $S$ is finite. 
    Since no variety in $S$ contains $\ab{A}$, all varieties
     in $S$ are strictly smaller than $\Var (\ab{A})$.
     We now show that for every variety $U$ with $U \subset \Var (\ab{A})$
     there is a variety $W \in S$ such that $U \subseteq W$.
     Since $U \subset \Var(\ab{A})$, \eqref{eq:WVB} yields $k$-ary terms $s, t$ such that
     $U \models s \approx t$ and $s^{\ab{A}} \neq t^{\ab{A}}$.
     Let $s_1, t_1 \in R$ with $s_1^{\ab{A}} = s^{\ab{A}}$ and $t_1^{\ab{A}} = t^{\ab{A}}$.
     Then, since $U \subseteq \Var (\ab{A})$, we have
     $U \models s_1 \approx s$ and $U \models t_1 \approx t$. Therefore
     $U \models s_1 \approx t_1$, and thus 
     $U \subseteq \Var (\ab{A}) \cap \Mod (s_1 \approx t_1)$, and this last variety 
     lies in $S$.

      Let $U$ be a variety with $U \lcover \Var (\ab{A})$.
     Then there is $W \in S$ with $U \subseteq W \subset \Var (\ab{A})$,
     and therefore $U = W$ and thus $U \in S$.
     Therefore all varieties $U$ with $U \lcover \Var(\ab{A})$ are among the
     finitely many elements of $S$; in fact, the cardinality
      of $S$ satisfies $|S| \le |\mathrm{Clo}_{k} (\ab{A})|^2$, where $\mathrm{Clo}_k(\ab{A})$ denotes
 the set of $k$-ary term functions of $\ab{A}$. This completes the proof
     of item~\eqref{it:co1}. 

     For proving~\eqref{it:co2}, we 
     let $U$ be a variety with $U \subset \Var(\ab{A})$. Let
     $W$ be maximal in $\{W' \in S \setsuchthat U \subseteq W'\}$. 
     We next show $W \lcover \Var (\ab{A})$. Suppose there is
     $W_1$ with $W \subset W_1 \subset \Var (\ab{A})$. Then there is $W_2 \in S$ with
     $W_1 \subseteq W_2 \subset \Var (\ab{A})$. Then $W \subset W_2$, contradicting the 
     maximality of $W$.
     Therefore $W \lcover \Var (\ab{A})$, and $U \subseteq W \lcover \Var (\ab{A})$. \qed

  \emph{Proof of Corollary~\ref{cor:4}:}
   \eqref{it:c1}$\Rightarrow$\eqref{it:c2}:
   Let $S := \{ W \setsuchthat W \le \Var (\ab{A}) \}$ be the class of subvarieties
   of $\Var (\ab{A})$.
   By  Theorem~\ref{thm:accsv2}, every variety in $S$ is finitely generated, and 
   so the assumptions yield that every variety in $S$ is finitely based relative
   to $\Var (\ab{A})$.
   From Lemma~\ref{lem:fb}, we obtain that $S$ has no infinite strictly
   descending chains. Thus, supposing that $S$ is infinite, we may pick
   $W$ minimal in $\{ W' \setsuchthat W'$ has infinitely many 
   subvarieties$\}$.
   Since $W$ is a subvariety of $\Var (\ab{A})$, it is finitely
   generated by Theorem~\ref{thm:accsv2}, and therefore, by Lemma~\ref{lem:finitelybranching},
   there must be a $W_1 \lcover W$ such that $W_1$ has infinitely
   many subvarieties. This contradicts the minimality of $W$.
    Hence $\Var (\ab{A})$ has only finitely many subvarieties.

  \eqref{it:c2}$\Rightarrow$\eqref{it:c3}:
    Let $(\ab{A}_i)_{i \in \N}$ be an infinite sequence of nonisomorphic cardinality
    critical algebras in $\Var (\ab{A})$.
    By the assumption, there is an infinite subset $T$ of $\N$ such that
    $\Var (\ab{A}_s) = \Var (\ab{A}_t)$ for all $s, t \in T$.
    We will now show that $|\ab{A}_s| = |\ab{A}_t|$ for all $s, t \in T$.
    Suppose $|\ab{A}_s| < |\ab{A}_t|$. 
 Since $\ab{A}_s \in \Var (\ab{A}_t)$, we have
    $\Var (\ab{A}_s) \subseteq \Var (\{ \ab{B} \setsuchthat \vb{B} \in \Var (\ab{A}_t), |\ab{B}| < | \ab{A}_t | \})$.
    Since $\ab{A}_t \in \Var (\ab{A}_s)$, this yields that $\ab{A}_t$ lies in the variety
      generated by the  members of $\Var (\ab{A}_t)$ of cardinality
    $< |\ab{A}_t|$, contradicting the fact that $\ab{A}_t$ is cardinality
    critical.
    Since $\Var (\ab{A})$ is locally finite, it contains only finitely many
    nonisomorphic members of the same finite cardinality because each $k$-element
    algebra is a homomorphic image of the free algebra $\ab{F}_{\Var (\ab{A})} (k)$.
    Hence $T$ must be finite, a contradiction.

 \eqref{it:c3}$\Rightarrow$\eqref{it:c4}:
     This follows from Lemma~\ref{lem:cc}.

 \eqref{it:c4}$\Rightarrow$\eqref{it:c1}:
     This follows from Lemma~\ref{lem:fb}. 

\qed
       

\begin{thebibliography}{OMVL78}

\bibitem[Aic10]{Ai:CMCO}
E.~Aichinger, \emph{Constantive {M}al'cev clones on finite sets are finitely
  related}, Proc. Amer. Math. Soc. \textbf{138} (2010), no.~10, 3501--3507.

\bibitem[AMM11]{AMM:OTNO}
E.~Aichinger, P.~Mayr, and R.~McKenzie, \emph{On the number of finite algebraic
  structures}, to appear in J. Eur. Math. Soc. (JEMS); available on arXiv:1103.2265v2 [math.RA].

\bibitem[BIM{\etalchar{+}}10]{BI:VWFS}
J.~Berman, P.~Idziak, P.~Markovi\'{c}, R.~McKenzie, M.~Valeriote, and
  R.~Willard, \emph{Varieties with few subalgebras of powers}, Transactions of
  the American Mathematical Society \textbf{362} (2010), no.~3, 1445--1473.

\bibitem[Bir35]{Bi:Comp}
G.~Birkhoff, \emph{On the structure of abstract algebras}, Proc. Cambridge
  Phil. Soc. \textbf{31} (1935), 433--454.

\bibitem[Bry82]{Br:TLOF}
R.~M. Bryant, \emph{The laws of finite pointed groups}, Bull. London Math. Soc.
  \textbf{14} (1982), no.~2, 119--123.

\bibitem[BS81]{BS:ACIU}
S.~Burris and H.~P. Sankappanavar, \emph{A course in universal algebra},
  Springer New York Heidelberg Berlin, 1981.

\bibitem[EN53]{EN:VGL}
T.~Evans and B.~H. Neumann, \emph{On varieties of groupoids and loops}, J.
  London Math. Soc. \textbf{28} (1953), 342--350.

\bibitem[Kru73]{Kr:ISFR}
R.~L. Kruse, \emph{Identities satisfied by a finite ring}, J. Algebra
  \textbf{26} (1973), 298--318.

\bibitem[Lee06]{Le:SVGF}
E.~W.~H. Lee, \emph{Subvarieties of the variety generated by the five-element
  {B}randt semigroup}, Internat. J. Algebra Comput. \textbf{16} (2006), no.~2,
  417--441.

\bibitem[{L'v}74]{Lv:VARI}
I.V. {L'vov}, \emph{{Varieties of associative rings. I.}}, {Algebra and Logic}
  \textbf{12} (1974), 150--167 (English).

\bibitem[McN92]{Mc:AFGT}
G.~F. McNulty, \emph{A field guide to equational logic}, J. Symbolic Comput.
  \textbf{14} (1992), no.~4, 371--397.

\bibitem[{Mil}85]{Mi:BWAB2}
E.C. {Milner}, \emph{{Basic wqo- and bqo-theory.}}, {Graphs and order. The role
  of graphs in the theory of ordered sets and its applications, Proc. NATO Adv.
  Study Inst., Banff/Can. 1984, NATO ASI Ser., Ser. C 147, 487-502},
  1985.

\bibitem[MMT87]{MMT:ALVV}
R.~N. McKenzie, G.~F. McNulty, and W.~F. Taylor, \emph{Algebras, lattices,
  varieties, volume {I}}, Wadsworth \& Brooks/Cole Advanced Books \& Software,
  Monterey, California, 1987.

\bibitem[Neu67]{Ne:VOG}
H.~Neumann, \emph{Varieties of groups}, Springer-Verlag New York, Inc., New
  York, 1967.

\bibitem[Ol{\cprime}70]{Ol:FBP}
A.~Ju. Ol{\cprime}{\v{s}}anski{\u\i}, \emph{The finite basis problem for
  identities in groups}, Izv. Akad. Nauk SSSR Ser. Mat. \textbf{34} (1970),
  376--384.

\bibitem[OMVL78]{OV:VTMO}
S.~Oates~MacDonald and M.~R. Vaughan-Lee, \emph{Varieties that make one
  {C}ross}, J. Austral. Math. Soc. Ser. A \textbf{26} (1978), no.~3, 368--382.

\bibitem[OP64]{OP:IRIF}
S.~Oates and M.~B. Powell, \emph{Identical relations in finite groups}, J.
  Algebra \textbf{1} (1964), 11--39.

\bibitem[PK79]{PK:FUR}
R.~P{\"o}schel and L.~A. Kalu{\v{z}}nin, \emph{Funktionen- und
  {R}elationenalgebren}, Mathematische Monographien, vol.~15, VEB Deutscher
  Verlag der Wissenschaften, Berlin, 1979.

\end{thebibliography}
\newcommand{\etalchar}[1]{$^{#1}$}
\def\cprime{$'$}
\providecommand{\bysame}{\leavevmode\hbox to3em{\hrulefill}\thinspace}
\providecommand{\MR}{\relax\ifhmode\unskip\space\fi MR }
\providecommand{\MRhref}[2]{%
  \href{http://www.ams.org/mathscinet-getitem?mr=#1}{#2}
}
\providecommand{\href}[2]{#2}

\begin{flushleft}
\begin{small}
Institut f\"ur Algebra, Johannes Kepler Universit\"at Linz,
Altenbergerstra\ss e 69,
4040
Linz,
Austria \\
{\tt erhard@algebra.uni-linz.ac.at} \\
{\tt peter.mayr@jku.at}\\[5mm]
\end{small}
\end{flushleft}
\end{document}